\documentclass[11pt]{article}
\usepackage{amssymb}
\usepackage{amsmath}

\newtheorem{theo}{Theorem}[section]

\newtheorem{prop}[theo]{Proposition}
\newtheorem{lemm}[theo]{Lemma}
\newtheorem{coro}[theo]{Corollary}
\newtheorem{rema}[theo]{Remark}
\newtheorem{Defi}[theo]{Definition}
\newtheorem{ex}[theo]{Example}

\newcommand{\cqfd}
{%
\mbox{}%
\nolinebreak%
\hfill%
\rule{2mm}{2mm}%
\medbreak%
\par%
}
\newfont{\gothic}{eufb10}
\date{\empty}
\begin{document}
\title{Hodge loci and absolute Hodge classes}
\author{Claire Voisin\\ Institut de math{\'e}matiques de Jussieu, CNRS,UMR
7586} \maketitle \setcounter{section}{-1}
\section{Introduction}
\setcounter{equation}{0}
Let $\pi:\mathcal{X}\rightarrow T$ be a
family of smooth projective complex varieties. Assume $\mathcal{X},
\pi, T$ are defined over $\mathbb{Q}$. An immediate consequence of
the fact that there are only countably many components of the
relative Hilbert scheme for $\pi$, and that the relative Hilbert
scheme (with fixed Hilbert polynomial) is defined over $\mathbb{Q}$,
is the following: if the Hodge conjecture is true, the components of
the Hodge locus in $T$ are defined over $\overline{\mathbb{Q}}$, and
their Galois transforms are again components of the Hodge locus. (We
recall later on the definition of the components of the Hodge
locus.) In \cite{CDK}, it is proven that the components of the Hodge
locus (and even the components of the locus of Hodge classes, which
is a stronger notion) are algebraic sets, while Hodge theory would
give them only a local structure of  closed analytic subsets (see
\cite{voisinbook}, 5.3.1).

In this paper, we give  simple sufficient conditions for components
of the Hodge locus to be defined over $\overline{\mathbb{Q}}$ (and
their Galois transforms to be also components of the Hodge locus).
This criterion of course does not hold in full generality, and it
particular does not say anything about the definition field of an
isolated point in the Hodge locus. But in practice, it is reasonably
easy to check and allows to conclude in some explicit cases, where
the Hodge conjecture is not known to hold. We give a few examples of
applications in section \ref{sec3}.

We will first relate this geometric language to the notion of
absolute Hodge classes (as we only deal with the de Rham version, we
will not  use the terminology of Hodge cycles of \cite{delmil}), and
explain why this notion allows to reduce the Hodge conjecture to the
case of varieties defined over $\overline{\mathbb{Q}}$, thus
clarifying a question asked to us by V. Maillot and Ch. Soul\'e.

Let us recall the notion of (de Rham) absolute Hodge class (cf
\cite{delmil}). Let $X^{an}$ be a complex projective manifold and
$\alpha\in Hdg^{2k}(X^{an})$ be a rational Hodge class. Thus
$\alpha$ is rational and

\begin{eqnarray}\label{comparaison30mai}
\alpha\in F^kH^{2k}(X^{an},\mathbb{C})\cong
\mathbb{H}^{2k}(X^{an},\Omega_{X^{an}}^{\bullet\geq k}).
\end{eqnarray}
Here, the left hand side is Betti cohomology of the complex manifold
$X^{an}$ and  the isomorphism of
(\ref{comparaison30mai}) is induced
by
the resolution
$$0\rightarrow {\mathbb C}\stackrel{(2i\pi)^k}{\rightarrow}\mathcal{O}
\stackrel{d}{\rightarrow}\Omega_X\rightarrow...\rightarrow\Omega_X^n\rightarrow 0,\,n=dim\,X$$
of the constant sheaf ${\mathbb C}$ on $X^{an}$.
The right hand side in (\ref{comparaison30mai}) can be computed, by GAGA
principle, as the hypercohomology of the algebraic variety $X$ with
value in the complex of algebraic differentials:
$$
\mathbb{H}^{2k}(X^{an},\Omega_{X^{an}}^{\bullet\geq k})\cong
\mathbb{H}^{2k}(X,\Omega_{X}^{\bullet\geq k}).$$
 Let us denote by
$\mathcal{E}$ the set of fields  embeddings of
$\mathbb{C}$ in $\mathbb{C}$. For each element
$\sigma $ of $\mathcal{E}$, we get a new algebraic
variety $X_\sigma$ defined over $\mathbb{C}$, and we have a similar
isomorphism for $X_\sigma$. Thus the class $\alpha$ provides a (de
Rham or Betti) complex cohomology class
$$\alpha_\sigma\in \mathbb{H}^{2k}(X_\sigma,\Omega_{X_\sigma}^{\bullet\geq
k})=F^k H^{2k}(X_\sigma^{an},\mathbb{C})
$$
for each  $\sigma\in \mathcal{E}$.
\begin{Defi} (cf \cite{delmil}) {\rm The class $\alpha $ is said to be (de Rham) absolute Hodge
if $\alpha_\sigma$ is a rational cohomology class for each
$\sigma$.}
\end{Defi}

We will introduce in section \ref{1} the notion of weakly absolute
Hodge class. In the definition above, we ask that each $\alpha_\sigma$ is proportional
to a rational cohomology class.

 We first prove in this Note
 the following statement,  which answers a question asked by Vincent Maillot and Christophe Soul\'e:
 \begin{prop} \label{prop}Assume the Hodge conjecture is known for varieties
 $X_{\overline{\mathbb{Q}}}$ defined over
 $\overline{\mathbb{Q}}$ and (weakly) absolute  Hodge classes
 $\alpha$ on them.
Then the Hodge conjecture is true for  (weakly) absolute  Hodge classes.
 \end{prop}
 \begin{rema}{\rm It is easy to see, (see
 Lemma \ref{leweak}) that a weakly absolute Hodge class $\alpha$ on a variety defined over
 $X_{\overline{\mathbb{Q}}}$ is
 defined over $\overline{\mathbb{Q}}$, that is
 $\alpha\in \mathbb{H}^{2k}(X_{\overline{\mathbb{Q}}},\Omega_{X_{\overline{\mathbb{Q}}}}^{\bullet\geq k})$.
 }
 \end{rema}
 \begin{rema} {\rm In the statement of the Proposition, we fix an embedding
 of $\overline{\mathbb{Q}}$ into $\mathbb{C}$, and so $\alpha$
 determines a class in
 $\mathbb{H}^{2k}(X_{\mathbb{C}},\Omega_{X_\sigma}^{\bullet\geq k})=F^kH^{2k}(X_{\mathbb{C}}^{an},\mathbb{C})$,
 which is assumed to be
 rational. If the Hodge conjecture
 is true for this class, then for any other embedding $\sigma$ of
 $\overline{\mathbb{Q}}$ into $\mathbb{C}$,  the class
 $\alpha_\sigma$ is also rational, and the Hodge conjecture is also true for this Hodge class.
 Thus the statement makes sense and is independent of the choice of embedding.}
 \end{rema}

 We next turn to the problem of whether Hodge classes should be absolute.
Let $X$ be a complex projective manifold, with a deformation family
  $\pi:\mathcal{X}\rightarrow T$  defined over $\mathbb{Q}$, that is
  $X$ is a fiber $X_0$, for some complex point $0\in T(\mathbb{C})$,
  and let $\alpha\in
 Hdg^{2k}(X)$ be a primitive  Hodge class.
We show the following:
 \begin{theo} \label{thm}1) Assume
 that for one irreducible component $S$ passing through $\alpha$  of the  locus of Hodge
 classes,
there is no constant sub-variation of Hodge structure of
$R^{2k}\pi_{S*}\mathbb{Q}_{prim} $ on $S$,
  except for $\mathbb{Q}\alpha_t$.
  Then $\alpha$ is weakly absolute.

 2) Let us weaken the assumptions on $S$ by asking that any
  constant sub-variation of
 Hodge structure of $R^{2k}\pi_{S*}\mathbb{Q}_{prim} $  on $S$  is of type $(k,k)$.
 Then, $p(S_{red})$ is defined over
 $\overline{\mathbb{Q}}$, and satisfies the property that its Galois translates are also
 of the form $p(S'_{red})$ for some irreducible  component $S'$ of the locus of Hodge classes.

 \end{theo}
 Here $\pi_S:\mathcal{X}_{S}\rightarrow S_{red}$ is obtained by base change
 $p: S_{red}\rightarrow T$.
 In statement 2), the Hodge locus of $\alpha $ is defined as the projection to
$T$ (via the projection map $p$) of the connected  component of the
locus of Hodge classes passing through $\alpha$. We will describe in
section \ref{1} their natural schematic structure.

Statement 1) will imply, by Lemma \ref{leweak} proven it next
section, that under the same assumptions,
  the  Hodge locus of $\alpha$ is defined over $\overline{\mathbb{Q}}$
  and its  image  under any element of
 $Gal\,(\overline{\mathbb{Q}}/\mathbb{Q})$ is again a component of the Hodge locus.

An immediate Corollary of Theorem \ref{thm},1) is the following
simple statement :
\begin{coro} Assume that the  infinitesimal Torelli theorem holds
for the variation of Hodge structure on
$R^{2k}\pi_{*}\mathbb{Q}_{prim} $. Assume
 that one component $S$ passing through $\alpha$  of the  locus of Hodge classes
  has positive dimension, and that the only proper non trivial sub-Hodge
  structure of $H^{2k}(X,\mathbb{Q})_{prim}$ is $\mathbb{Q}\alpha$.
  Then $\alpha$ is weakly absolute.
\end{coro}
Note that the assumption that $S$ has positive dimension is
satisfied once $h^{k-1,k+1}:=rk\,H^{k-1,k+1}(X)_{prim}<dim\,T$ (cf
\cite{voisinbook}, Proposition 5.14).

\vspace{0.5cm}

{\bf Proof.} Indeed, a constant sub-variation of Hodge structure of
$R^{2k}\pi_{S*}\mathbb{Q}_{prim} $  on $S$ must then be (by taking
the fiber at the point $0$ corresponding to $X$) either equal to
$R^{2k}\pi_{S*}\mathbb{Q}_{prim} $ or to $\mathbb{Q}\alpha$. The
first case is impossible by the Torelli assumption, and $dim\,S>0$.
Thus the assumptions of Theorem \ref{thm},1) are satisfied. \cqfd
Case 2) of Theorem \ref{thm}  leads to the following generalization of Proposition \ref{prop}:
\begin{prop}\label{prop'} Suppose the Hodge conjecture is true for Hodge classes
on smooth projective varieties defined over $\overline{\mathbb{Q}}$.
Then under the assumptions of Theorem \ref{thm}, 2), the class $\alpha$ is algebraic.
\end{prop}

Section \ref{1} is devoted to the discussion of absolute and weakly absolute
Hodge classes in terms of the corresponding components of the locus of Hodge classes
and components of the Hodge locus.

In section \ref{sec2}, we prove the results stated in this introduction.

We  give in the last section variants and  applications of Theorem
\ref{thm}. In Theorem \ref{th3}, we give an algebraic (Zariski open)
criterion on a Hodge class $\alpha\in F^kH^{2k}$ in order that the
assumptions of Theorem \ref{thm} are satisfied at least at a general
point of the connected component $\widetilde{S}_\alpha$ of the locus
of Hodge classes passing through $\alpha$. Of course, except in
level $2$, where we can use the Green density criterion, it is hard
to decide if there are many Hodge classes in the Zariski open set of
$F^kH^{2k}$ where this criterion is satisfied. We give examples in
level $2$, where this criterion is satisfied in a Zariski dense open
set, in which there are ``many'' Hodge classes. In one of these
examples, the Hodge conjecture is not known to hold for these
classes.

The second application (Theorem \ref{hypersurfaces}) concerns the period map. Under a reasonable
assumption on the infinitesimal variation of Hodge structures on the
primitive cohomology of the fibers of a family $\pi:\mathcal{X}\rightarrow T$
of projective varieties defined over $\mathbb{Q}$, we conclude that
any component $W$ dominating $T$ by the first projection of the set
of pairs $(t,t')\in T\times T$ such that the Hodge structures on
$H^n(X_t,\mathbb{Q})_{prim}$ and $H^n(X_{t'},\mathbb{Q})_{prim}$ are isomorphic,
is defined over $\overline{\mathbb{Q}}$.
\section{Absolute and weakly absolute Hodge classes\label{1}}
\setcounter{equation}{0}
Let us introduce the following  variant of
the notion of absolute Hodge class.
\begin{Defi} {\rm The class $\alpha $ is said to be weakly (de Rham) absolute Hodge
if for each  $\sigma\in \mathcal{E}$, $\alpha_\sigma$ is a multiple
$\lambda_\sigma\gamma_\sigma$, where  $\gamma_\sigma\in
H^{2k}(X^{an},\mathbb{Q})$ is a rational cohomology class (hence a
Hodge class) and $\lambda_\sigma\in \overline{\mathbb{Q}}$.}
\end{Defi}
\begin{rema}\label{rem}{\rm It turns out that the condition
$\lambda_\sigma\in \overline{\mathbb{Q}}$ is automatically satisfied. Indeed, consider
the primitive decomposition of $\alpha$ with respect to the polarization
given by the projective embedding of $X$.
$$\alpha=\sum_{2k-2r\geq 0,2r\leq n}c_1(L)^{k-r}\alpha_r,\,n=dim\,X,$$
where $\alpha_r\in H^{2r}(X,\mathbb{Q})_{prim}$.

Then the primitive decomposition of $\alpha_\sigma$ is given by
$$\alpha_\sigma=\sum_{2k-2r\geq 0,2r\leq n}c_1(L_\sigma)^{k-r}\alpha_{r,\sigma},$$
and thus, if
$\alpha_\sigma=\lambda_\sigma\gamma_\sigma$, with $\gamma_\sigma\in
H^{2k}(X_\sigma^{an},\mathbb{Q})$, then for each $r$, we get
\begin{eqnarray}\label{eqn1}\alpha_{r,\sigma}=\lambda_\sigma\gamma_{\sigma,r},
\end{eqnarray}
where $\gamma_{\sigma,r}$ is the degree $2r$ primitive component of
$\gamma_\sigma$, and thus is a rational cohomology class. But we
know by the second Hodge-Riemann bilinear relations that if
$\alpha_r\not=0$, we have $\int_Xc_1(L)^{n-2r}\alpha_r^2\not=0$.
This is a rational number, which is also equal to
$\int_{X_\sigma}c_1(L_\sigma)^{n-2r}\alpha_{r,\sigma}^2$. On the
other hand, as $\gamma_{\sigma,r}$ is a rational cohomology class,
we also have
$\int_{X_\sigma}c_1(L_\sigma)^{n-2r}\gamma_{\sigma,r}^2\in
\mathbb{Q}$, and thus, from the equalities \begin{eqnarray}
\label{eqn3}
\int_{X}c_1(L)^{n-2r}\alpha_{r}^2=\int_{X_\sigma}c_1(L_\sigma)^{n-2r}\alpha_{r,\sigma}^2=\lambda_\sigma^2
\int_{X_\sigma}c_1(L_\sigma)^{n-2r}\gamma_{\sigma,r}^2,
\end{eqnarray} we get $\lambda_\sigma^2\in \mathbb{Q}$. }
\end{rema}
Geometrically, the meaning of these notions is the following
 (see also \cite{esnaultparanjape}):
Let $\pi:\mathcal{X}\rightarrow T$ be a family of deformations of
$X$, which is defined over $\mathbb{Q}$ (here $T$ is not supposed to
be geometrically irreducible, and thus the assumption is not
restrictive on $X$). There is then the algebraic  vector bundle
$F^kH^{2k}$ on $T$,  defined over $\mathbb{Q}$, which is the total
space of the locally free sheaf
$F^k\mathcal{H}^{2k}=R^{2k}\pi_*\Omega_{\mathcal{X}/T}^{\bullet\geq
k}$ on $T$. We will use the following terminology (see \cite{CDK}).
The {\it locus of Hodge classes} for the  family above, in degree
$2k$
 is the
   set
 of pairs $(X_t,\alpha_t)\in F^kH^{2k}(\mathbb{C})$,
such that $\alpha_t\in H^{2k}(X_t^{an},\mathbb{C})$ is rational (hence a Hodge class).

 The components of the  {\it Hodge locus}  are  the image
  in $T$, via the natural projection $p:F^k{H}^{2k}\rightarrow T$, of the connected components
  of the  locus of Hodge classes. If $\alpha$ is a Hodge class on $X$, the
  {\it Hodge locus of} $\alpha$ is the image in $T$ of the connected component
of the  locus of Hodge classes passing through $\alpha$.

  Notice that the  locus of Hodge classes is obviously locally a
  countable union of closed analytic subsets in $F^kH^{2k}(\mathbb{C})$
   Indeed, if $\alpha\in F^kH^{2k}(X^{an}_t,\mathbb{C})\cap H^{2k}(X^{an},\mathbb{Q})$,
  then in a simply connected neighbourhood $U$ of $t\in T$, we
  have a trivialization of the locally constant sheaf
  $R^{2k}\pi_*^{an}\mathbb{C}$,  which induces a trivialization of the corresponding
  vector bundle $H^{2k}$ and  gives a composed holomorphic map:
  \begin{eqnarray}
  \label{psi}\psi: F^kH^{2k}\hookrightarrow {H}^{2k}\rightarrow
  H^{2k}(X_t,\mathbb{C}),
  \end{eqnarray}
  where
  ${H}^{2k}$ is the
total space of the locally free sheaf
$$\mathcal{H}^{2k}=R^{2k}\pi_*\Omega_{\mathcal{X}/T}^{\bullet}=R^{2k}\pi_*^{an}\mathbb{C}\otimes\mathcal{O}_T$$ on
$T$.

  Then, over $U$, the locus of Hodge classes identifies to
  $\psi^{-1}(H^{2k}(X_t,\mathbb{Q}))$, which is a countable union of fibers of $\psi$.
  This defines a natural schematic structure on the connected components of the locus of Hodge classes.

  Similarly, the local description of the Hodge locus of $\alpha$ is as follows: we
   can locally extend
  $\alpha$ to a locally constant section $\tilde{\alpha}$ of $R^{2k}\pi^{an}_*\mathbb{Q}$.
  Then $\tilde{\alpha}$ gives in particular a holomorphic section
  of the vector bundle $\mathcal{H}^{2k}:=R^{2k}\pi^{an}_*\mathbb{C}\otimes \mathcal{O}_T$.
  Then the Hodge locus of $\alpha$ is simply defined by the condition
  $$\tilde{\alpha}^{\leq k-1}=0,$$
  where $\tilde{\alpha}^{\leq k-1}$ is the projection of $\tilde{\alpha}$ in the quotient
  $\mathcal{H}^{2k}/F^k\mathcal{H}^{2k}$. This again defines the
  schematic structure of the Hodge locus of $\alpha$.

It is  clear from these descriptions that the projection from the  locus of Hodge classes to the
Hodge locus is a local immersion which is open onto an union of  local analytic branch of the Hodge locus.

    Cattani, Deligne and Kaplan proved in fact the following much stronger result
    concerning the structure of the locus of Hodge classes
    (cf \cite{CDK}):
    \begin{theo}\label{thcdk} The connected components of the locus of Hodge classes
    are algebraic subsets of the algebraic vector bundle $F^kH^{2k}$.
    \end{theo}

  The conjecture that Hodge classes are absolute Hodge is equivalent to saying that
  the  locus of Hodge classes is a countable union of algebraic subsets defined over $\mathbb{Q}$. To see
  this, note that
  a given Hodge class is absolute Hodge if and only if its $\mathbb{Q}$-Zariski closure
  in $F^kH^{2k}$
  is contained in the locus of Hodge classes. It is then clear by the Noetherian
  property  that for countably many
  generically chosen Hodge classes $\alpha_i$, the  locus of Hodge classes must be equal
  to the union of the $\mathbb{Q}$-Zariski closures of the $\alpha_i$'s.

  The statement that Hodge classes are weakly absolute Hodge implies
   the facts that the locus of Hodge classes is a countable union
  of algebraic subsets of $F^kH^{2k}$ defined over $\overline{\mathbb{Q}}$
  and that
 the  Hodge locus is a countable union of algebraic subsets
  of $T$ defined over $
  {\mathbb{Q}}$.  More precisely, we have :
  \begin{lemm} Let $\alpha\in H^{2k}(X^{an},\mathbb{Q})$ be a weakly absolute Hodge
  class. Then the connected component $\widetilde{S}_\alpha$
  of the locus of Hodge classes
  passing through $\alpha$ is defined (schematically) over
  $\overline{\mathbb{Q}}$, and so is the Hodge locus of $\alpha$.
  Furthermore the Galois images of the Hodge locus of $\alpha$ are
  also (schematically) components of the Hodge locus.\label{leweak}
  \end{lemm}
{\bf Proof. } We know by Theorem \ref{thcdk} that $\widetilde{S}_{\alpha}$ is
algebraic, and it is by definition connected. We make the base
change $\widetilde{S}_{\alpha,red}\rightarrow T$, where
we replace if necessary $\widetilde{S}_{\alpha}$ by a Zariski open set, in order
to make the reduced scheme $\widetilde{S}_{\alpha,red}$ smooth. Then the corresponding
family
$$\pi_\alpha:\mathcal{X}_\alpha\rightarrow \widetilde{S}_{\alpha,red}$$
admits the locally constant section $\tilde\alpha\in
H^0(\widetilde{S}_{\alpha,red},R^{2k}\pi_{\alpha*}\mathbb{Q})$. By the
global invariant cycle theorem \cite{de2}, there exists a class
$\beta\in H^{2k}(\overline{\mathcal{X}}_\alpha,\mathbb{Q})$ which is
of type $(k,k)$ and restricts to $\tilde\alpha_t$ on  each fiber $X_t$ of
the family $\mathcal{X}_\alpha$. In fact, we can even make this
class canonically defined by choosing an ample line bundle
$\mathcal{L}$ on $\overline{\mathcal{X}}_\alpha$, which allows to
define a polarization on $
H^{2k}(\overline{\mathcal{X}}_\alpha,\mathbb{Q})$ (see also the
proof of Proposition \ref{prop} or \cite{andre} for more details).
Then $\beta$ is canonically defined if we impose that $\beta$ lies
in the orthogonal complement of  $Ker\,rest_X$ with respect to this
polarization.

Let now $\alpha$ be weakly absolute. Then the class $\beta_\sigma$
on $\mathcal{X}_{\alpha,\sigma}$ restricts to
$\alpha_\sigma=\lambda_\sigma\gamma_\sigma$ on $X_\sigma$, where
$\gamma_\sigma$ is rational, and is in the orthogonal complement of
$Ker\,rest_{X_\sigma}$ with respect to the polarization induced by
$\mathcal{L}_\sigma$.
  It thus follows that
$\frac{1}{\lambda_\sigma}\beta_\sigma$, which restricts to
$\gamma_\sigma$, has to be rational (hence is a Hodge class).
Let $\tilde{\gamma}$ be the
locally constant section of
$R^{2k}\pi_{\alpha*}\mathbb{Q}$ on $\widetilde{S}_{\alpha,red}$ obtained by restricting
$\frac{1}{\lambda_\sigma}\beta_\sigma$.
  We conclude that we have an
  inclusion

\begin{eqnarray}
\label{eqn2}\frac{1}{\lambda_\sigma}\sigma(\widetilde{S}_{\alpha,red})\subset\widetilde{S}_{\gamma_\sigma,red}
,\end{eqnarray} which is easily checked to extend in fact to a schematic
identification
\begin{eqnarray}
\label{eqn2222}\frac{1}{\lambda_\sigma}\sigma(\widetilde{S}_{\alpha})=\widetilde{S}_{\gamma_\sigma}
.\end{eqnarray}
Indeed, this follows  from the flatness of the sections $\tilde{\alpha}_\sigma$ and
$\tilde\gamma_\sigma$,   from the fact that
$\lambda_\sigma$ has to be constant  along
$\sigma(\widetilde{S}_{\alpha,red})$ by formula (\ref{eqn3}), and from the fact that
$\widetilde{S}_{\gamma_\sigma}$ is by definition connected.

As
$$p(\frac{1}{\lambda}\sigma(\widetilde{S}_\alpha))=p(\sigma(\widetilde{S}_\alpha)))=
\sigma(p(\widetilde{S}_\alpha)),$$
we conclude from (\ref{eqn2222}) that the image via $\sigma$ of the
Hodge locus of $\alpha$ is also a component of the Hodge locus.

Finally, to see that if $\alpha$ is weakly absolute Hodge, then
$\widetilde{S}_\alpha\subset  F^kH^{2k}$ is defined over $\overline{\mathbb{Q}}$,
  we
  use equality
  (\ref{eqn2222}), applied to
  $\sigma\in \mathcal{E}$ together with the fact noticed in Remark
  \ref{rem} that $\lambda^2_\sigma\in\mathbb{Q}$. It follows that
  the constant
  $\lambda_\sigma\in\overline{\mathbb{Q}}$
  can take only countably many values, and
  in particular, there are only countably many Galois transforms
 $\sigma(\widetilde{S}_\alpha )$, and as we know that
$\widetilde{S}_\alpha$ is algebraic, this implies that
$\widetilde{S}_\alpha$ is defined over $\overline{\mathbb{Q}}$.

\cqfd

\section{Proof of Theorem \ref{thm} and Propositions \ref{prop}, \ref{prop'}.\label{sec2}}
\setcounter{equation}{0}

{\bf Proof of Proposition \ref{prop}.} Let $(X^{an},\,\alpha)$ be a
pair consisting of a projective complex manifold and an
 absolute (resp. a weakly absolute)
rational Hodge class. By the geometric interpretation given above,
and by Lemma \ref{leweak} in the weakly absolute case, it follows
that there exist smooth irreducible quasi-projective varieties
$\mathcal{X},\,T$ defined over $\overline{\mathbb{Q}}$, a projective
morphism $\pi:\mathcal{X}\rightarrow T$, and a locally constant
global section
$$\tilde{\alpha}\in H^0(T,R^{2k}\pi_*\mathbb{Q}),$$
such that $X$ is one fiber of $\pi$ and $\alpha$ is the restriction of
$\tilde{\alpha}$ to this fiber.

   Deligne's global invariant cycle theorem  \cite{de2} says now that for any smooth
    completion
    $\overline{\mathcal{X}}$ of $\mathcal{X}$, there exists a Hodge class
    $\beta \in Hdg^{2k}(\overline{\mathcal{X}})$ such that
    $$\beta_{\mid X}=\alpha.$$

    Of course, we may also choose $\overline{\mathcal{X}}$ defined over $\overline{\mathbb{Q}}$.
In order to conclude, we claim that we may choose $\beta$
to be absolute Hodge (resp. weakly absolute Hodge). Indeed, we will
deduce from this, under the assumptions of Proposition \ref{prop},
that
 $\beta$ is the class of an algebraic cycle, and then, so is its restriction $\alpha$.

To prove the claim, consider the morphism of rational
Hodge structures
$$H^{2k}(\overline{\mathcal{X}}^{an},\mathbb{Q})\rightarrow H^{2k}(X^{an},\mathbb{Q}).
$$
The left hand side can be polarized using a ample line bundle
$\mathcal{L}$ on $\overline{\mathcal{X}}$. (That is, we use the Lefschetz
decomposition with respect to this polarization, and change the
signs of the natural intersection pairing
$$(\alpha_r,\beta_r)=\int_{\overline{\mathcal{X}}}c_1(\mathcal{L})^{N-2r}\alpha_r\cup
\beta_r,\,N=dim\,\mathcal{X}$$ on the pieces of the Lefschetz
decomposition with $r$ even, in order to get a polarized Hodge
structure.) Thus we conclude that there is an orthogonal
decomposition
$$H^{2k}(\overline{\mathcal{X}}^{an},\mathbb{Q})=A\oplus B$$
into the sum of two Hodge structures, where the first one identifies
via restriction to its image in $H^{2k}(X,\mathbb{Q})$, while the
second one is the kernel of the restriction map.
$B$ is a sub-Hodge
structure of $H^{2k}(\overline{\mathcal{X}}^{an},\mathbb{Q})$
and
$A$ is then defined as the orthogonal of $B$ under the metric
described above on $H^{2k}(\overline{\mathcal{X}}^{an},\mathbb{Q})$.

We define then $\beta$ to be the unique element of $A$ which
restricts to $\alpha$.

For each element
$\sigma$ of $Gal\,(\overline{\mathbb{Q}}/\mathbb{Q})$, we get a line bundle
$\mathcal{L}_\sigma$ on
$\overline{\mathcal{X}}^{an}_\sigma$,  a sub-Hodge structure
$B_\sigma:=Ker\,rest_{X_\sigma}$, and the isomorphism
\begin{eqnarray} H^{2k}(\overline{\mathcal{X}}^{an},\mathbb{C})\cong
H^{2k}(\overline{\mathcal{X}}^{an}_\sigma,\mathbb{C})
\label{compar28}
\end{eqnarray}
commutes with restrictions maps and is compatible with the
polarizations given by
$\mathcal{L}$ and $\mathcal{L}_\sigma$.
Thus we get similarly  a rational sub-Hodge structure
$A_\sigma$ of  $H^{2k}(\overline{\mathcal{X}}^{an}_\sigma,\mathbb{Q})$
and there is a commutative diagram where the horizontal maps are restrictions maps
 and
thus are defined on rational cohomology, and the vertical maps are induced  by the comparison
isomorphism (\ref{compar28}) :
$$\begin{matrix}A\otimes\mathbb{C}&\hookrightarrow &H^{2k}(X,\mathbb{C})\\
\shortparallel&&\shortparallel\\
A_\sigma\otimes\mathbb{C}&\hookrightarrow &H^{2k}(X_\sigma,\mathbb{C})
\end{matrix}.$$
It follows from this
that if $\alpha$ is absolute Hodge (resp. weakly absolute Hodge), so
is $\beta$.

\cqfd

{\bf Proof of Theorem \ref{thm}.} 1) Let $(X,\alpha)$ be as in the
statement of the Theorem. By Theorem \ref{thcdk},  the component
passing through $\alpha$ of the locus of Hodge classes is an
algebraic set. Let $S$ be an irreducible component of this set
containing $(X,\alpha)$, and satisfying the assumption of Theorem
\ref{thm},1). Replacing $S$ by a Zariski open set of $S_{red}$, we
may assume  that $S$ is smooth. There is by base change a projective
family $\pi_S:\mathcal{X}_S\rightarrow S$ together with a tautological flat
section
$$\tilde{\alpha}\in H^0(S,R^{2k}\pi_{S*}\Omega_{\mathcal{X}_S/S}^{\bullet\geq k}),$$
with value $\alpha_t$ at each $t$.

Let $\overline{\mathcal{X}}_S$ be a smooth completion of
$\mathcal{X}_S$. The global invariant cycle theorem says that there
exists a class $\beta \in
H^{2k}(\overline{\mathcal{X}}_S,\mathbb{Q})\cap
F^kH^{2k}(\overline{\mathcal{X}}_S,\mathbb{Q})$ such that
$\beta_{\mid X}=\alpha$. On the other hand, the vector space
$$ H^{2k}(\overline{\mathcal{X}}_S,\mathbb{Q})_{\mid X_t}\cap H^{2k}(X_t,\mathbb{Q})_{prim}$$
is a constant sub-Hodge structure of $H^{2k}(X,\mathbb{Q})_{prim}$.
Thus, by our assumption on $S$, we conclude that it must be equal to
$\mathbb{Q}\alpha_t$. It follows that the complex vector space
$$ \mathbb{H}^{2k}(\overline{\mathcal{X}}_S,\Omega_{\overline{\mathcal{X}}_S}^\bullet)_{\mid X_t}\cap
\mathbb{H}^{2k}(X_t,\Omega_{X_t}^\bullet)_{prim}$$ has rank $1$ and
is generated by $\alpha$.

 Let $\sigma\in \mathcal{E}$. We want
to show that  the class
$$\alpha_\sigma\in
\mathbb{H}^{2k}(X_\sigma,\Omega_{X_\sigma}^{\bullet\geq k})\subset
H^{2k}(X_\sigma^{an},\mathbb{C})$$ is of the form
$\lambda_\sigma\gamma_\sigma$, where $\gamma_\sigma$ is rational.

But $\sigma $ provides a new family
$\overline{\mathcal{X}}_{S,\sigma}$ fibered over $S_\sigma$ with
fiber $X_{t,\sigma}$, such that the vector space
\begin{eqnarray}\label{restDR}
\mathbb{H}^{2k}(\overline{\mathcal{X}}_{S,\sigma},\Omega_{\overline{\mathcal{X}}_{S,\sigma}}^\bullet)_{\mid
X_\sigma} \cap
\mathbb{H}^{2k}(X_\sigma,\Omega_{X_\sigma}^\bullet)_{prim}
\end{eqnarray}
has rank $1$ and is generated by $\alpha_\sigma$. It follows that
the intersection of the image of the restriction map
\begin{eqnarray}\label{restbetti} H^{2k}(\overline{\mathcal{X}}_{S,\sigma}^{an},\mathbb{Q})\rightarrow
H^{2k}(X_\sigma^{an},\mathbb{Q}),
\end{eqnarray}
 with $H^{2k}(X_\sigma^{an},\mathbb{Q})_{prim}$ has rank
$1$.

Thus we have $\alpha_\sigma=\lambda_\sigma\gamma_\sigma$ for some
rational primitive Hodge class $\gamma_\sigma$ on $X_\sigma$, and
some non zero complex coefficient $\lambda_\sigma$. By Remark
\ref{rem}, we have $\lambda_\sigma\in \overline{\mathbb{Q}}$, and
thus $\alpha$ is weakly absolute.

2) The proof of 2) is very similar. Indeed, with the  same notations as above,
we find that $\alpha$ belongs to the sub-Hodge structure
$$H^{2k}(\overline{\mathcal{X}}_S^{an},\mathbb{Q})_{\mid X^{an}}\cap
H^{2k}(X^{an},\mathbb{Q})_{prim},$$ which is the fiber at $0$ of the
locally constant
 sub-Hodge structure
$$H^{2k}(\overline{\mathcal{X}}_S^{an},\mathbb{Q})_{\mid X_t^{an}}\cap
H^{2k}(X_t^{an},\mathbb{Q})_{prim},\,t\in S,$$ hence must be a
trivial sub-Hodge structure. This assumption is algebraic, as it can
be translated into the fact that the vector space
$$\mathbb{H}^{2k}(\overline{\mathcal{X}}_S,\Omega_{\overline{\mathcal{X}}_S}^\bullet)_{\mid X}\cap
\mathbb{H}^{2k}(X,\Omega_{X}^{\bullet})_{prim}$$ is equal to
$$\mathbb{H}^{2k}(\overline{\mathcal{X}}_S,\Omega_{\overline{\mathcal{X}}_S}^{\bullet\geq k})_{\mid X}\cap
\mathbb{H}^{2k}(X,\Omega_{X}^{\bullet\geq k})_{prim}.$$

Let now $\sigma\in \mathcal{E}$. We conclude from
the above  that the sub-Hodge structure
$$H^{2k}(\overline{\mathcal{X}}_{S,\sigma}^{an},\mathbb{Q})_{\mid X_\sigma^{an}}\cap
H^{2k}(X_\sigma^{an},\mathbb{Q})_{prim},
 $$
 to which $\alpha_\sigma$ belongs,
 is trivial.
 Thus we can write
$\alpha_\sigma=\sum_{i=1}^{i=N}\alpha_i\gamma_i$, where $\gamma_i$
are independent rational Hodge classes on $X_\sigma^{an}$ coming
from $\overline{\mathcal{X}}_{S,\sigma}^{an}$
 and the $\lambda_i$ are complex coefficients.
As $\alpha_\sigma$ gives a flat section of
 the bundle $F^k\mathcal{H}^{2k}$ on $\sigma(S)$,
and the $\gamma_i$ are locally constant  on $\sigma(S_{red})$, we
conclude that
 the $\lambda_i$'s are constant on $\sigma(S_{red})$. We claim
 that for a generic choice of {\it rational}
 coefficients $\lambda'_i,\,1\leq i\leq N$, $\sigma(p(S_{red}))$
 is equal to
 $p(S''_{red})$ where
 $S''$ is an irreducible component of  the  locus of Hodge classes passing through
 $\sum_{i=1}^{i=N}\lambda'_i\gamma_i$.

 Assuming the claim, this shows that there are only countably many
 Galois transforms of $p(S_{red})$ and thus, because $p(S_{red})$ is
 algebraic, this implies that $p(S_{red})$ is defined over
 $\overline{\mathbb{Q}}$. The claim also gives the second part of
 the statement.

 To prove the claim, we choose a simply connected neighborhood $U$ of
 the point $\sigma(0)\in T(\mathbb{C})$.
 Over $U$, we can consider the
 map $\psi:F^kH^{2k}_{\mid U}\rightarrow H^{2k}(X_{\sigma}^{an},\mathbb{C})$ of
 (\ref{psi}). Then for any choice of complex coefficients $\mu_i$,
 we know that $p(\psi^{-1}(\sum_i\mu_i\gamma_i))$ contains
 $p(\sigma(S_{red}))\cap U$, and that for $(\mu_1,\ldots,\mu_N)=(\lambda_1,\ldots,\lambda_N)$,
  $p(\sigma(S_{red}))\cap U$  is the reduction of an irreducible component
  of $p(\psi^{-1}(\sum_i\mu_i\gamma_i))$. By lower semi-continuity of
  the dimension of the fibers of $\psi$, we conclude that the later property remains true
  for
  $(\mu_i)\in\mathbb{C}^N$ in a Zariski open set of coefficients,
  and thus in particular for some $N$-uple
  $(\lambda'_i)\in\mathbb{Q}^N$.

Having this, we proved that for some irreducible analytic component
$S'$ of $\psi^{-1}(\sum_{i=1}^{i=N}\lambda'_i\gamma_i)$, the two
analytic subsets  $\sigma(p(S_{red}))\cap U$ and $p(S'_{red})$ of $U$ coincide.
Because $\sigma(p(S_{red}))$ is irreducible and reduced, and  because
by Theorem \ref{thcdk}, $\psi^{-1}(\sum_{i=1}^{i=N}\lambda'_i\gamma_i)$ is an
open set in an irreducible algebraic subset $S''$ of $F^kH^{2k}$, (an
irreducible component of a connected component of the locus of Hodge classes), we get by
analytic continuation that
  $\sigma(p(S_{red}))=p(S''_{red})$.

\cqfd
\begin{rema} The schematic structure of the  locus where  a
combination $\sum_i\mu_i\gamma_i$ remains in $F^kH^{2k}$ may depend
on the $\mu_i$, even if we know that the corresponding reduced
algebraic set does not depend generically on the $\mu_i$'s. This is why we have to restrict here
to the reduced subschemes.
\end{rema}
Let us conclude this section by giving the proof of
proposition \ref{prop'}.

\vspace{0,5cm}

{\bf Proof of Proposition \ref{prop'}.} Indeed, with the same notations as above,
we just proved that $S':=p(S_{red})$ is defined over
$\overline{\mathbb{Q}}$. We also know that the only locally constant sub-Hodge structure
of $R^{2k}\pi_*\mathbb{Q}_{prim}$ is of type $(k,k)$.
As monodromy acts in a finite way on
the set of  Hodge classes of $X_t,\,t\in S'$ generic,
there is an \'etale cover $S''$ of the smooth part of $S'$, also defined
over $\overline{\mathbb{Q}}$, on which this monodromy action becomes trivial.
Thus we have by base change a family
$\pi'':\mathcal{X}_{S''}\rightarrow S''$, together with a
 global section $\tilde{\alpha}$
of $R^{2k}\pi''_*\mathbb{Q}_{prim}$, whose restriction to $X_0$ is equal to $\alpha$.
The global invariant cycle theorem now says that
there exists a Hodge class $\beta$ on a smooth compactification
$ \overline{\mathcal{X}}_{S''}$, which we may assume  defined over $\overline{\mathbb{Q}}$,
restricting to $\alpha$.
If the Hodge conjecture is true for Hodge classes on varieties defined over
$\overline{\mathbb{Q}}$, it is then true for $\beta$ and thus also for $\alpha$.
\cqfd

\section{\label{sec3} Variants and applications}
\setcounter{equation}{0}

Let us give to start with an infinitesimal criterion which will
guarantee that the assumptions of Theorem \ref{thm}, 1) are
satisfied by an irreducible component of  $\widetilde{S}_\alpha$.
 This will then give as a consequence of Theorem \ref{thm} an
algebraic criterion (Theorem \ref{th3}) for a Hodge class $\alpha\in
F^kH^{2k}(\mathbb{C})$, to be
 weakly absolute.

 We assume again that $\pi:\mathcal{X}\rightarrow T$ is  a family of
 projective varieties defined over $\mathbb{Q}$, and we denote
 as before by $F^kH^{2k}$ the algebraic vector bundle whose sheaf of
 sections is
 $(R^{2k}\pi_*\Omega_{\mathcal{X}/T}^{\bullet\geq k})_{prim}$, which
 admits as a quotient the bundle $H^{k,k}$ whose sheaf of sections
 is $(R^{k}\pi_*\Omega_{\mathcal{X}/T}^{
 k})_{prim}$. This is an algebraic vector bundle defined over
 $\mathbb{Q}$.
 We have the $\mathcal{O}_T$-linear map which describes the infinitesimal variation of Hodge structure
 $$\overline{\nabla}:\mathcal{H}^{p,q}\rightarrow
 \mathcal{H}^{p-1,q+1}\otimes\Omega_T,$$
 which is defined using the Gauss-Manin connection and Griffiths
 transversality (cf \cite{voisinbook}, 5.1.2).
 Here $\mathcal{H}^{p,q}:=(R^{q}\pi_*\Omega_{\mathcal{X}/T}^{
 p})_{prim}$.

 The assumption of positive dimension for the Hodge loci
  is automatically satisfied if $h^{k-1,k+1}(X)_{prim}<dim\,T$. This is proved
 in \cite{voisinbook}, where it is shown that the Hodge loci in $T$ for the variation of Hodge
 structure on  $H^{k-1,k+1}(X_t)_{prim}$ can be defined by at most
 $h^{k-1,k+1}(X)_{prim}<dim\,T$.
We assume below that $T$ is smooth.

 Let $\alpha\in H^{2k}(X,\mathbb{Q})_{prim}$ be a Hodge
 class, where $X=X_0$ is a fiber of $\pi$, $0\in T(\mathbb{C})$.
 Let $\lambda\in H^{k,k}$ be the projection of
 $\alpha\in F^kH^{2k}$ in $H^{k,k}$.

   Let us assume that  the
  map
  $$\mu:T_{T,0}\rightarrow H^{k-1,k+1}(X_0)$$
  given by $\mu(v)=\overline{\nabla}_v(\lambda)$ is surjective. Let
  $K_\lambda$ be its kernel. $K_\lambda$ is the tangent space of
  the Hodge locus of $\alpha$ at $0$ (cf \cite{voisinbook}???).
We have the following
 algebraic  criterion on $\lambda$, for $\alpha$ to be weakly
 absolute:
 \begin{theo}\label{th3} Assume that

 1) $\mu_\lambda$ is surjective.

 2)  For
 $p>k,\,p+q=2k$ the map
 $$\overline{\nabla}_0:H^{p,q}(X_0)_{prim}\rightarrow H^{p-1,q+1}(X_0)\otimes
 K_\lambda^*,$$
 obtained by restriction of $\overline{\nabla}$, is injective.

 3) The map $$H^{k,k}(X_0)_{prim}\rightarrow H^{k-1,k+1}(X_0)\otimes
 K_\lambda^*,$$
 obtained by restriction of $\overline{\nabla}$, has for kernel the line generated by
 $\lambda$.

 Then $\alpha$ is weakly absolute.
 \end{theo}
{\bf Proof.} As the map $\mu$ is surjective,
the component $S_\alpha$ of the Hodge locus determined by $\alpha$
is smooth with tangent space $K_\lambda$ at $0\in T$ (cf
\cite{voisinbook}, Proposition 5.14).

The conditions 2) and 3) imply that any constant sub-variation of
Hodge structure of $R^{2k}\pi_*\mathbb{Q}_{prim}$  defined along an
open set of $S_\alpha$ containing the point $0$ parameterizing $X$
is equal to $\mathbb{Q}\alpha$. Indeed, if $\gamma^{p,q}$ is a
locally constant section of $R^{2k}\pi_*\mathbb{C}_{prim}$ which
remains of type $(p,q)$ near $0$ on ${S}_\alpha$, where we may
assume $p\geq q$ by Hodge symmetry, then we have
$$\overline{\nabla}\gamma^{p,q}=0\,{\rm in}\,
H^{p-1,q+1}\otimes\Omega_{S_\alpha},$$ and in particular, we have at
$0$, $$\overline{\nabla}\gamma^{p,q}(0)=0\,{\rm in}\,
 H^{p-1,q+1}(X_0)_{prim} \otimes K_\lambda^*.$$
 Thus by assumptions 2) and 3), we conclude that
 $\gamma^{p,q}=0$ for $p>k$ and $\gamma^{p,q}$ is proportional to
 $\lambda$ for $p=k$.

 We conclude then by applying  Theorem \ref{thm}, 1).
 \cqfd
\begin{rema} The same reasoning shows that if we only assume 1) and 2) in Theorem
\ref{th3},  then the class $\alpha$ satisfies the conclusion of
part 2 of Theorem \ref{thm}.  Thus in particular, if $\widetilde{S}_\alpha^0$ is the irreducible
component of $\widetilde{S}_\alpha$ passing through $\alpha$ (it is unique and reduced because
$\widetilde{S}_\alpha$ is now smooth at the point $\alpha$), then
$p(\widetilde{S}_\alpha^0)$ is defined over $\overline{\mathbb{Q}}$.
\end{rema}

 It is interesting  to note that condition 1) is a Zariski open
 condition  on the class
 $\lambda\in H^{k,k}$ (non necessarily Hodge)
and that conditions 2) and 3) are Zariski open in the set where
 1) is sastified. One can even note that the complementary set where
 these conditions are not satisfied, is Zariski closed and defined
 over $\mathbb{Q}$ as are the bundles $\mathcal{H}^{p,q}$ and the
 map $\overline{\nabla}$.

 Of course, even if we can show that the Zariski open set of
 $F^kH^{2k}$ defined by the condition 1), 2), 3) above is non empty,
 it is not clear if there are any Hodge classes in it. This is the
 case however if our variation of Hodge structure has Hodge level
 $2$, that is $h^{p,q}=0$ for $p\geq k+2$. Indeed, in this case, we
 have the Green density criterion (cf \cite{voisinbook}, 5.3.4) which
 guarantees that if there is any $\lambda\in H^{k,k}(X_t)_{prim}$
 satisfying property 1), then the set of rational Hodge classes are
 topologically dense in the real part of the vector bundle
 $H^{k,k}$.

\begin{ex} The criterion above allows us to prove that many Hodge classes
are
 weakly absolute for surfaces in $\mathbb{P}^3$, without using
 the Lefschetz theorem on $(1,1)$-classes.
 \end{ex}
 More interestingly, it allows to show a similar result for
 certain level $2$ subvariations of Hodge structure in the
 $H^{2k}$ of a variety $X$, without knowing the Hodge-Grothendieck
 generalized conjecture for this sub-Hodge structure.
 We can construct such examples on $4$-dimensional hypersurfaces
 with automorphisms.
 \begin{ex} Consider the action
 of the involution
 $\iota$ on $\mathbb{P}^5$ given by
 $\iota(X_0,\ldots,X_5)=(-X_0,-X_1,X_2,\ldots,X_5)$, and
 take for $T$ the family of isomorphism classes of degree $6$ hypersurfaces whose defining
 equation is invariant under $\iota$, and for
 sub-Hodge structure the anti-invariant part of $H^4(X)$ under
 $\iota$. This Hodge structure has Hodge level $2$, because $\iota$
 acts trivially on the rank $1$ space $H^{4,0}(X)$.
 Thus the Green
 density criterion applies once assumption 1) above is satisfied.
  The parameter space
 $T$ has dimension $226 $, and the number $h^{1,3}_-$ is equal to $208 $.
One can check that assumptions 1), 2), 3) are satisfied generically
on
  $F^2H^4_-$, thus proving
 that many Hodge classes are weakly absolute, even if the Hodge
 conjecture is not known for them. This is done following \cite{cagr}, \cite{grirat} by computations
 in the Jacobian ring of the generic hypersurface as above. In fact this can be done for
 $X$ the Fermat hypersurface, and for a generic class
 $\lambda\in H^{2,2}(X)_-$.

 \end{ex}

We now turn to another application of Theorem \ref{thm}, which
concerns  the fibers of the period map and the Torelli problem.

Let $\pi:\mathcal{X}\rightarrow T$ be a family of smooth projective
varieties which is defined over $\mathbb{Q}$, and consider the
variation of Hodge structure on $H^n(X_t)_{prim}$. Here $T$ is
assumed to be smooth. The corresponding infinitesimal variation of
Hodge structure at $t\in T$ is given by the map
$$\overline{\nabla}: H^{p,q}(X_t)_{prim}\rightarrow
H^{p-1,q+1}(X_t)_{prim}\otimes \Omega_{T,t}.$$
 We will
assume the following property : at the generic point $0\in T$, the corresponding map
$$\mu: H^{p,q}(X_0)_{prim}\otimes T_{T,0}\rightarrow
H^{p-1,q+1}(X_0)_{prim},$$ $$\eta\otimes v\mapsto
\overline{\nabla}_v(\eta),$$
 is surjective whenever
$H^{p,q}(X_0)_{prim}\not=0$. (This property is satisfied for example
by the families of hypersurfaces or complete intersections in
projective space.) We then have :
\begin{theo} \label{hypersurfaces} Let $Z \subset T\times T$ be
the set of points $(t,t')$ such that there exists an isomorphism of
Hodge structures between $H^{n}(X_t,{\mathbb Q})_{prim}$ and $
H^{n}(X_{t'},\mathbb{Q})_{prim}$.
 Let $W\subset T\times T$ be (the underlying reduced scheme
 of) an irreductible component of $Z$ which dominates $T$.
 Then under the assumptions above, $W$ is defined over $\overline{\mathbb{Q}} $ and  any Galois
 transform of $W$ is again (the underlying reduced scheme
 of) an irreducible component of $Z$.
\end{theo}
{\bf Proof.} We apply Theorem \ref{thm}, 2). The set $W$ above is
$\Gamma_{red}$ for an irreducible component $\Gamma$ of the Hodge
locus corresponding to the induced variation of Hodge structure of weight $0$ on
$H^{n}(X_t,\mathbb{Q})_{prim}^*\otimes
H^{n}(X_{t'},\mathbb{Q})_{prim}$ on  $T\times T$. What we have to
prove in order to apply Theorem \ref{thm},2) is the fact  that if
$W$ dominates $T$ by the first (or equivalently second) projection,
then any constant sub-Hodge structure of
$$H^{n}(X_t,\mathbb{Q})_{prim}^*\otimes
H^{n}(X_{t'},\mathbb{Q})_{prim},\,\,(t,t')\in W,$$
 must be of type $(0,0)$.

By definition, for $(t,t')\in W$, the Hodge structures on
$$H^{n}(X_t,\mathbb{Q})_{prim},\,\,\,
H^{n}(X_{t'},\mathbb{Q})_{prim},$$ are isomorphic. Thus the Hodge
structures on $$ H^{n}(X_t,\mathbb{Q})_{prim}^*\otimes
H^{n}(X_{t'},\mathbb{Q})_{prim},\,\,\,
H^{n}(X_t,\mathbb{Q})_{prim}^*\otimes H^{n}(X_{t},\mathbb{Q})_{prim}$$
are isomorphic.
Furthermore, $t$ is generic. Thus it suffices to
prove that on any finite cover of $T$,  there is no constant sub-Hodge structure of
$H^{n}(X_t,\mathbb{Q})_{prim}^*\otimes H^{n}(X_{t},\mathbb{Q})_{prim}$
which is not of type $(0,0)$.

This is done by an easy infinitesimal argument. Let $$
\alpha\in
H^{n}(X_t,\mathbb{Q})_{prim}^*\otimes
H^{n}(X_{t},\mathbb{Q})_{prim}$$ be of bidegree $(r,s)$ with $r>s,\, r+s=0$. Thus
$r>0$,  and if we see
$\alpha$ as an element of $Hom\,(H^{n}(X_t)_{prim},H^{n}(X_t)_{prim})$,
 $\alpha\in H^{r,-r}$ means that
$\alpha(H^{p,q}(X_t)_{prim})\subset H^{p+r,q-r}(X_t)_{prim}$.

We have to show that if there is a flat section
 $\tilde{\alpha}$ on $T$, extending $\alpha$ and staying of type
 $(r,-r)$, then $\alpha=0$.
 It suffices to show this at first order at $0\in T$, where this is equivalent
 to say that if
 $$\overline{\nabla}\alpha=0\in H^{r-1,-r+1}(X_0\times
 X_0)\otimes\Omega_{T,0},$$
 then $\alpha=0$.
 But to say that $\overline{\nabla}\alpha=0$ is equivalent to say that
 \begin{eqnarray}\label{eqn14}
 \overline{\nabla}_v(\alpha(\phi))=\alpha(\overline{\nabla}_v(\phi)),\,\forall \phi\in
 H^{p,q}(X_0)_{prim},\,\forall (p,q),\,p+q=n,\,\forall v\in T_{T,0}.
 \end{eqnarray}
  Equation (\ref{eqn14}) shows
that $\alpha$ is in fact determined by its value on the first non
$0$ term $H^{p,q}(X_0)_{prim}$, because by assumption the map
$$H^{p,q}(X_0)_{prim}\otimes T_{T,0}\rightarrow
H^{p-1,q+1}(X_0)_{prim}, $$ $$ \phi\otimes v\mapsto
\overline{\nabla}_v(\phi)$$ is surjective once $H^{p,q}(X_0)_{prim}$
is different from $0$.

On the other hand, $\alpha$ must be zero on the first non $0$ term
$H^{p,q}(X_0)_{prim}$, because it sends it in
$H^{p+r,q-r}(X_0)_{prim}=0$.
 \cqfd

\end{document}